\documentclass[reqno,a4]{amsart}
\usepackage{amsmath,amssymb,amsthm}
\usepackage{graphicx}

\newcommand{\R}{{\mathbb R}}
\newcommand{\Q}{{\mathbb Q}}
\newcommand{\Z}{{\mathbb Z}}
\newcommand{\T}{{\mathbb T}}
\newcommand{\C}{{\mathbb C}}
\newcommand{\N}{{\mathbb N}}
\newcommand{\A}{{\mathcal A}}
\newcommand{\B}{{\mathcal B}}
\newtheorem{theorem}{Theorem}
\newtheorem{lem}[theorem]{Lemma}
\newtheorem{cor}[theorem]{Corollary}

\theoremstyle{definition}
\newtheorem{ex}[theorem]{Example}
\newtheorem{rem}[theorem]{Remark}

\newcommand{\I}{\mathfrak{i}}
\newcommand{\e}{\mathbf{e}}
\renewcommand{\d}{\mathbf{d}}

\begin{document}
\title{Asymptotic formula for balanced words}
\author{Shigeki Akiyama}
\maketitle
\begin{abstract}
We give asymptotic formulas for the number of balanced words whose slope $\alpha$ and intercept $\rho$
lie in a prescribed rectangle. They are related to uniform distribution of Farey fractions and Riemann Hypothesis.
In the general case, the error term is deduced using an inequality of large sieve type.
\end{abstract}

\section{Introduction}
Let $\A=\{0,1\}$ and we denote by $\A^*$ the monoid generated by $\A$ by concatenation,
where the empty word $\lambda$ is its identity. The length of $x\in \A^*$ is $|x|$. 
We denote $|x|_1$ the cardinality of $1$ in $x\in \A^*$.
$\A^{\N}$ (resp. $\A^{\Z}$) is the set of right infinite  (resp. bi-infinite) words. 
If $y\in \A^*$ is a factor (a subword) of $x\in \A^{*}\cup \A^{\N} \cup \A^{\Z}$, we write $y\prec x$.
A word in $x\in \A^*\cup \A^{\N} \cup \A^{\Z}$ is {\bf balanced} if $||u|_1-|v|_1|\le 1$ holds 
for any $u,v\prec x$ with $|u|=|v|$.
An infinite word $w\in \A^{\N}$ is {\bf sturmian} if $\text{Card} \{ u\prec w\ |\ |u|=n\}=n+1$ for all $n\in \N$. 
Morse and Hedlund characterized a sturmian word as an aperiodic balanced word in $\A^{\N}$. They
also characterized sturmian words as a coding of irrational rotation. More precisely
a lower mechanical word $(s_n)\in \A^{\N}$  is defined by
$$
s_n(\alpha, \rho)=
\lfloor \alpha(n+1)+\rho\rfloor -
\lfloor \alpha n+\rho\rfloor 
$$
with a given slope $\alpha\in [0,1]$ and an intercept $\rho\in [0,1)$. 
An upper mechanical word is similarly defined by replacing $\lfloor \cdot \rfloor$ with $\lceil \cdot \rceil$,
which is denoted by $\hat{s}_n(\alpha,\rho)$.
Then a sturmian word is a (lower or upper) mechanical word of an irrational slope $\alpha$ and vice versa. 
It is known that every balanced word is a factor of a lower and an upper mechanical word \cite{Lothaire:02}. 
Indeed, for every balanced word $x=x_1\dots x_n$ we can find
a slope $\alpha$ and an intercept $\rho$ such that $x_i=s_i(\alpha, \rho)=\hat{s}_i(\alpha,\rho)$.
Note that 
the choice of $\alpha$ and $\rho$ is not unique. In fact, a balanced $x$ word corresponds to $(\rho,\alpha)$
in a convex polygon in $[0,1)\times [0,1]$. This geometric idea for enumeration 
is found in \cite{Berstel-Pocchiola, Yasutomi:98, Kamae-Takahashi, Cassaigne-Frid}.
Let $\phi(n)$ be the Euler totient function. 
The formula for the number of balanced words of length $n$ is given by
\begin{equation}
\label{Balanced}
1+\sum_{k=1}^n (n+1-k)\phi(k).
\end{equation}
Several different proofs are found in 
\cite{Lothaire:02, Lipatov, Berstel-Pocchiola, ADLuca-Mignosi, Motomaki-Saari}.

In this paper we refine the formula (\ref{Balanced}) and give its asymptotic behavior.
Denote by $B(n,t,u)$ the cardinality of the set of balanced words of length $n$ whose slope $\alpha\in
[1-t,1]$\footnote{The statement of Theorem \ref{Count} is simpler by this choice than taking $[0,t]$
because the lines in the proof do not intersect $(0,1)\times \{0\}$.}
and its intercept $\rho\in [0,u)$ and let $\B(n)$ be the set of balanced words of length $n$.
Then we show
\begin{theorem}
\label{Count}
$$
B(n,t,u)=1+ \sum_{m\le n} A(m,t,u)
$$
with
$$
A(m,t,u)=\sum_{\substack{i<j\le m,\ (i,j)=1\\ i/j\le t,\ \langle mi/j \rangle < u}} 1,
$$
where $i$ and $j$ are non negative integers. 
\end{theorem}
Theorem \ref{Count} 
slightly generalizes Yasutomi \cite[Proposition 4]{Yasutomi:98}, shown in a different context.
By using Theorem \ref{Count}, we will derive asymptotic formulas.
Hereafter, we use conventional terminology in analytic number theory, i.e., Landau $O$, $o$
symbol and Vinogradov symbol $\ll$. 
The symbol $\varepsilon$ is reserved as an arbitrary positive constant, and the symbol
$c$ in Landau $O$ is a suitably chosen positive constant, which may differ among formulas.
\begin{theorem}
\label{Main}
$$
B(n,t,u)= \frac{tu}{\pi^2} n^3 + O\left(n^{2} (\log n)^{15/2}\right).
$$
Moreover, we have
$$
B(n,1,1)=\frac {n^3+3n^2}{\pi^2}+O\left(n^2\exp\left(-c\left((\log n)^{3/5}(\log \log n)^{-1/5}\right)\right)\right)
$$
and
\begin{equation}
\label{Gen}
B(n,t,1)=\frac {tn^3}{\pi^2}+O(n^2).
\end{equation}
\end{theorem}
For almost all $t$, the estimate (\ref{Gen}) can be sharpened with the help of 
Fujii \cite{Fujii:07}, see the discussion in the end of \S 3.
Denote by $\chi_A(x)$ the indicator function of the set $A$. 
For a rectangle $(a,b]\times [c,d)$ in the unit square $[0,1)\times [0,1]$, we see
$$
\mathrm{Card} \{ x\in \B(n)\ |\ (\rho,\alpha)\in (a,b]\times [c,d) \} = \frac{(b-a)(d-c)}{\pi^2} n^3+ O\left(n^{2}(\log n)^{15/2}\right),
$$
from
$$
\chi_{(a,b]\times [c,d)}(x)=
\chi_{(0,b]\times [c,1)}(x)-
\chi_{(0,a]\times [c,1)}(x)-
\chi_{(0,b]\times [d,1)}(x)+
\chi_{(0,a]\times [d,1)}(x).
$$
Moreover for a Jordan measurable region $W$ in the unit square, we have
\begin{cor}
$$
\mathrm{Card} \{ x\in \B(n)\ |\ (\rho,\alpha)\in W \} = \frac{\mathrm{Area}(W)}{\pi^2} n^3+ O\left(n^{2}(\log n)^{15/2}\right)
$$
\end{cor}
\noindent
where $\mathrm{Area}$ is the $2$-dimensional Lebesgue measure, 
since every Jordan measurable set is well approximated by finite union of rectangles.

Farey series $(f_m(i))$ of order $m$ is the finite increasing sequence composed of irreducible
 fractions in $[0,1)$ whose denominators are not larger than $m$:
$$
0=f_m(1)<f_m(2)<\dots <f_m(\Phi(m))<1
$$
with $\Phi(m)=\sum_{k=1}^m \phi(k)$. Clearly $A(m,t,1)=\max \left\{j\ |\ f_m(j) \le t \right\}$.
It is well-known that Riemann Hypothesis is equivalent to
$$
\sum_{i=1}^{\Phi(m)} \exp(2\pi\I f_m(i))\ll m^{1/2+\varepsilon},
$$
a strong uniform distribution property of Farey fractions.
As pointed out in \cite{Fujii:05}, Franel \cite{Franel} already noticed that
$$
\int_{0}^1 (A(m,t,1)-t \Phi(m))^2 dt \ll m^{1+\varepsilon}
$$
is also equivalent to Riemann Hypothesis. One can see
$$
A(m,t,1)-t \Phi(m)=O(m),
$$
similarly to (\ref{u1}) and (\ref{Est0}) in \S 3, but we expect it is much smaller in average. 
See also \cite{Niederreiter:73, Fujii:05, Kanemitsu-Yoshimoto:96}.
We state another equivalent statement directly related to the number of balanced words:
\begin{cor}
\label{RH}
The estimate
\begin{equation}
\label{RHB}
B(n,1,1)=\frac {n^3+3n^2}{\pi^2}+O\left(n^{3/2+\varepsilon}\right)
\end{equation}
is equivalent to Riemann Hypothesis.
\end{cor}

\section{Proof of Theorem \ref{Count}}
\begin{figure}[htb]
\includegraphics[width=0.7\columnwidth]{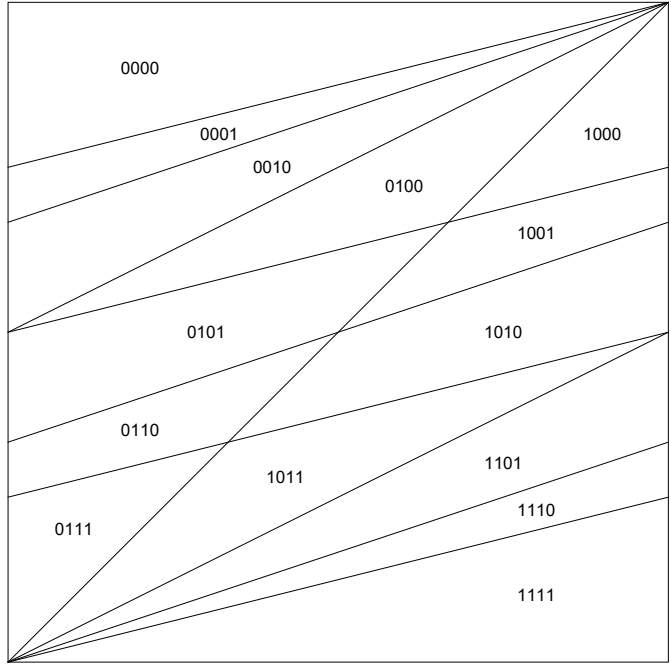}
\caption{Partition for $m=4$ \label{Dis}}
\end{figure}

We elucidate a
geometric counting discussion of Yasutomi \cite[Proposition 4]{Yasutomi:98} in our convenient
 terminology, which is more straightforward than the one in \cite{Berstel-Pocchiola}.
Let $m$ be a fixed positive integer and put $X:=[0,1)\times [0,1]$.
The map $\psi:(\rho,\alpha) \to (s_n(1-\alpha,\rho))_{n=1}^m$ gives
a natural partition
\begin{equation}
\label{Part}
X := \bigcup_{w\in B(m)} \psi^{-1}(w).
\end{equation}
Then $X$ is cut into convex cells by segments:
$$
Y:=X \setminus \{ (x,y) \ |\ x=n y-\ell, n\in \{1,\dots,m\}, \ell\in \{0,1,\dots, n-1\} \}.
$$
We obtain essentially the same partition by using $\hat{s}_n$,  the difference is seen 
only on the boundary of $Y$. In this paper, we use $s_n$ for the partition.
Fixing $\alpha\in [0,1]$, the intersections of 
the line $y=1-\alpha$ and $x=n y-\ell$ are
written as
$$
\{ R^{-n}(0) \ |\ n\in \{1,\dots,m\} \}\times \{1-\alpha\} =\{ -n\alpha\, (\bmod{1})\ |\ n=1,2\dots,m \}\times \{1-\alpha\}
$$
where $R:(x,1-\alpha)\mapsto (x+\alpha,1-\alpha)$ is the rotation map acting on the unit interval $[0,1)\times \{1-\alpha\}$ 
which is identified
with the torus $\T:=\R/\Z$.
The partition of $[0,1)\times \{1-\alpha\}$ 
by $R^{-n}(0)\ (n=1,\dots,m)$ gives cylinder sets of rotation $R$, i.e., the points in the 
same cylinder share the same coding 
$$\left(\chi_{[1-\alpha,1)\times \{1-\alpha\}}\left(R^{i-1}(x)\right)\right)_{i=1}^{m}.$$
They are the
$m+1$ different words which correspond to the words of length $m$ appear
in the sturmian word of slope $\alpha$ when $\alpha$ is irrational.
Slicing $X$ by the line $y=1-\alpha$, 
we observe the $m$-th level cylinder sets 
of the rotation $R$ on $\T$ of slope $\alpha$ acting on 
the unit interval $[0,1)\times \{1-\alpha\}$.
Considering $\alpha$ as a variable, we reconstruct the partition (\ref{Part}) of $X$ that every
convex cell corresponds to an element of $\B(m)$, 
which is consistent with the cylinder
partition of $[0,1)\times \{1-\alpha\}$ for each $\alpha\in [0,1]$.
In this manner, the partition (\ref{Part}) is seen as a pile of cylinder sets of 
level $m$ for all (rational/irrational) rotations.
The case $m=4$ is depicted in Figure \ref{Dis}.

To enumerate $B(n,t,u)$, we compute 
$B(m,t,u)-B(m-1,t,u)$ for $m\ge 2$, that is, the increase of number of cells as we add a new slope $1/m$. 
Equivalently, we count the number of intersections
in $[0,u)\times [0,t]$ which appear by adding new segments of slope $1/m$, see Figure \ref{Ex}.
The intersection points $(x,y)\in [0,u)\times [0,t]$ of $x=my-b$ and $x=\ell y-c\ (\ell<m)$ 
are in one to one correspondence with the set of their $y$-coordinates:
$$
P(m):=\left\{ \frac{b-c}{m-\ell}\ \left|\ \begin{aligned}  & \frac bm \le \frac {b-c}{m-\ell} < \frac {b+1}m,\quad
\frac{b-c}{m-\ell}\le t,\\ & \frac{\ell b-mc}{m-\ell}<u,\quad 0\le b,\ell<m,\quad 0\le c<\ell 
\end{aligned}
\right.\right\}.
$$
We claim that $P(m)$ coincides with the set
$$
Q(m):=\left\{ x\in \Q\cap [0,1) \ \left|\ x\le t,\quad \langle mx \rangle < u,\quad \mathrm{den}(x)\le m \right. \right\}
$$
where $\langle x\rangle=x-\lfloor x\rfloor$ and $\mathrm{den}(y)$ is the denominator of $y\in \Q$.
In fact $P(m) \subset Q(m)$ is clear. 
Let $x=p/q\in Q(m)$ with $1\le q\le m$, $(p,q)=1$, $p/q\le t$ and $\langle mx\rangle<u$.
Put $\ell=m-q$ and choose $b$ with $b/m\le p/q<(b+1)/m$ and let $c=b-p$. Since $p/q<(b+1)/m$, we have $p\le b$ and $c\ge 0$.
From the property of Farey fraction, we have $(b-p)/(m-q)< b/m < p/q$, which implies $c=b-p< m-q=\ell$. The inequality
$0\le \frac{\ell b-mc}{m-\ell}<1$ follows from $\frac bm \le \frac {b-c}{m-\ell} < \frac {b+1}m$. Thus 
$\langle mx \rangle= \frac{\ell b-mc}{m-\ell}$ and therefore $x\in P(m)$. The claim is proved. 
Set $A(m,t,u)=\mathrm{Card}(Q(m))$ for $m\ge 1$. From $B(1,t,u)=2$ and $A(1,t,u)=1$,
we obtain Theorem \ref{Count}.\qed
\medskip

We illustrate this proof by an

\begin{ex}
Let $t=0.7$ and $u=0.59$.  Then $(A(m,t,u))_{m=1}^{8}=(1, 2, 4, 4, 7, 8, 10, 13)$
and $B(8,t,u)=50$.  In Figure \ref{Ex}, 
we are counting the number of cells in the shaded region $[0,u)\times [0,t]$. 
Dashed segments are
of slope $1/8$ intersecting at 13 points 
indicated by dots, which contribute to the increase of the number of cells.

\begin{figure}[htb]
\includegraphics[width=0.7\columnwidth]{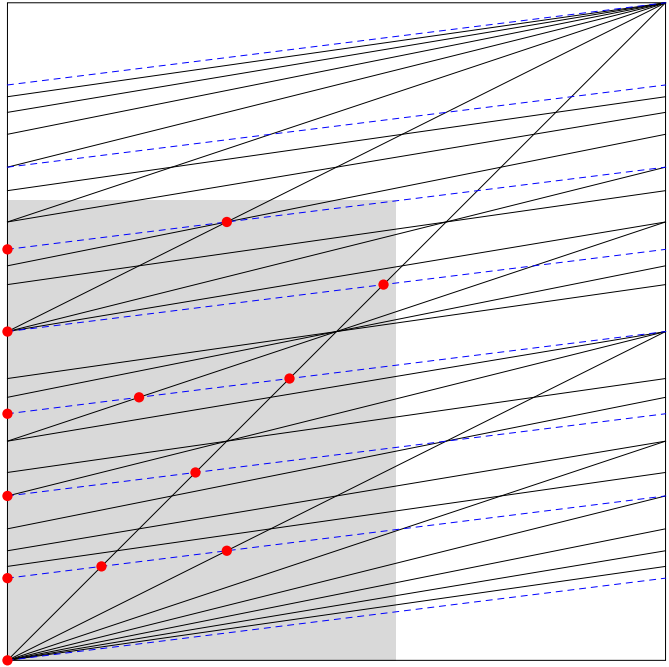}
\caption{$B(8,0.7, 0.59)=50$ \label{Ex}}
\end{figure}
\end{ex}

Hereafter we show an asymptotic formula of $B(n,t,u)$. 

\section{Theorem \ref{Main} for $u=1$ and Corollary \ref{RH}}
We show a
\begin{lem}
\label{Cesaro}
$$
\sum_{m\le x} (x-m)\phi(m)= \frac{x^3}{\pi^2} + O\left(x^2\exp\left(-c(\log x)^{3/5}(\log \log x)^{-1/5}\right)\right).
$$
\end{lem}

\proof
This is a variation of 
the argument to deduce prime number theorem \cite[Theorem 6.9]{Montgomery-Vaughan}.
We use the Mellin inversion formula
$$
\sum_{m\le x} (x-m)\phi(m)=\frac{1}{2\pi\I} \int_{a-\I \infty}^{a+\I \infty}
\frac {\zeta(s-1)x^{s+1}}{\zeta(s)s(s+1)}ds
$$
for $a>2$. Here $\zeta(s)$ is Riemann zeta function with a complex variable $s=\sigma+ \I w$
and put $\tau=|w|+4$. Since $|\zeta(s-1)|\ll \tau^{1/2}\log \tau$ and $1/\zeta(s)\ll \log \tau$
in the required range,
shifting the path to $\sigma=1+(\log x)^{-1}=:a_0$, we can pick the residue at $s=2$:
\begin{equation}
\label{Shift1}
\sum_{m\le x} (x-m)\phi(m)=\frac {x^3}{\pi^2}+
\frac{1}{2\pi\I} \int_{a_0-\I \infty}^{a_0+\I \infty}
\frac {\zeta(s-1)x^{s+1}}{\zeta(s)s(s+1)}ds.
\end{equation}
Truncating this formula, 
$$
\sum_{m\le x} (x-m)\phi(m)=\frac {x^3}{\pi^2}+\frac{1}{2\pi\I} \int_{a_0-\I T}^{a_0+\I T}
\frac {\zeta(s-1)x^{s+1}}{\zeta(s)s(s+1)}ds + O\left(\frac{x^{2}(\log T)^2}{T^{1/2}}\right)
$$
and shifting the path to the zero-free region 
$$\sigma>1-\frac{c}{(\log \tau)^{2/3}(\log\log \tau)^{1/3}}=:b$$
of $\zeta(s)$ due to
Vinogradov-Korobov \cite{Walfisz, Ivic},
we get a rectangular contour to be studied.
The contribution from horizontal segments is $O(x^{2}/T^{3/2-\varepsilon})$
and the vertical one at $\sigma=b$ is $O(x^{b+1})$.
Lemma \ref{Cesaro} is obtained by the choice $T=\exp(c (\log x)^{3/5}(\log\log x)^{-1/5})$.
\qed

Let us start with the easiest case $t=u=1$.
From 
$$
A(m,1,1)=\sum_{\substack{i<j\le m\\(i,j)=1}} 1=\Phi(m)
$$
and Theorem \ref{Count}, we have
\begin{eqnarray}
\label{Cum}
B(n,1,1)&=&1+\sum_{m=1}^n \Phi(m)\\
\nonumber
&=& 1+\sum_{j=1}^n (n+1-j)\phi(j)\\
\label{B11}
&=&
\frac {n^3+3n^2}{\pi^2} + O\left(n^2\exp\left(-c(\log n)^{3/5}(\log \log n)^{-1/5}\right)\right).
\end{eqnarray}
Here we used Lemma \ref{Cesaro} and the Mertens formula
$$
\Phi(n)= \frac {3n^2}{\pi^2}+ O(n \log n).
$$
The error term $E(n)=\Phi(n)-(3/\pi^2)n^2$ is well
studied in literature \cite{Pillai-Chowla,ErdosShapiro:51,Walfisz,Montgomery:87,Kaczorowski-Wiertelak}. 
However in (\ref{Cum}), the effect of this error term cancels out and we find
the second main term in (\ref{B11}). We retrieved the formula (\ref{Balanced}) as well.

Now we show Corollary \ref{RH}. If Riemann Hypothesis is valid, then (\ref{Shift1}) holds with 
$a_0=1/2+\varepsilon$. Thus from
$$
\lim_{T\to \infty}
\int_{a_0-\I T}^{a_0+\I T}\left| \frac{\zeta(s-1)}{\zeta(s)s(s+1)}\right|ds < \infty,
$$
we get the estimate (\ref{RHB}). 
Conversely Mellin transformation shows
$$
\int_{1}^{\infty} \left(\sum_{n\le x} (x-n)\phi(n)-\frac{x^3}{\pi^2}\right)
x^{-s-2} dx = \frac{\zeta(s-1)}{\zeta(s)s(s+1)} - \frac 1{\pi^2(s-2)}
$$
for $\sigma>2$. If (\ref{RHB}) is valid, then the parenthesis in the integrand is $O(x^{3/2+\varepsilon})$. 
This gives 
the holomorphic continuation of the right side to $\sigma>1/2+\varepsilon$, which finishes the proof.


Let us discuss a small counting issue. 
Since $i=0$ implies $j=1$ in the sum $A(m,t,u)$, we have
$$
A(m,t,u)=1+\sum_{\substack{i<j\le m,\, i/j\le t\\\langle mi/j \rangle < u,\, (i,j)=1} }1
$$
where $i$ and $j$ are {\bf positive} integers hereafter. 
Therefore
$$
B(n,t,u)=1+n+
\sum_{\substack{i<j\le m\le n,\, i/j\le t\\\langle mi/j \rangle < u,\, (i,j)=1} }1.
$$
If $t<1$ then, it is the same as 
$$
B(n,t,u)=1+n+
\sum_{\substack{i,j\le m\le n,\, i/j\le t\\\langle mi/j \rangle < u,\, (i,j)=1} }1.
$$
In the case $t=1$, $i=j$ happens only when $i=j=1$, and we may write 
$$
B(n,1,u)=1+\sum_{\substack{i,j\le m\le n\\\langle mi/j \rangle < u,\, (i,j)=1} }1.
$$
Let $\mu$ be the M\"obius function and assume $u=1$ and $t<1$. With positive integers $a, b$, we have
\begin{eqnarray}
\nonumber
B(n,t,1)-1-n&=&\sum_{\substack{i,j\le m\le n\\ i/j\le t}}\, \sum_{k \mid (i,j)} \mu(k)
=\sum_{k\le m\le n} \mu(k) \sum_{b\le m/k}\, \sum_{a/b\le t} 1\\
\nonumber
&=&\sum_{k\le m\le n}\mu(k) \sum_{b\le m/k} \lfloor bt \rfloor\\
\label{u1}
&=&t \sum_{k\le m\le n}\mu(k) \sum_{b\le m/k}b
-\sum_{kb\le m\le n}\mu(k)\langle bt \rangle.
\end{eqnarray}
While $t=1$, the same computation gives
$$
B(n,1,1)-1=\sum_{k\le m\le n}\mu(k) \sum_{b\le m/k}b,
$$
we see
\begin{equation}
\label{T1_0}
B(n,t,1)=1-t+n +t B(n,1,1)-\sum_{kb\le m\le n}\mu(k)\langle bt \rangle
\end{equation}
for $t<1$. 
Since 
\begin{equation}
\label{Est0}
\sum_{kb\le m}\mu(k)\langle bt \rangle=O(m)
\end{equation}
by
Niederreiter \cite[(6) and Lemma 3]{Niederreiter:73},
we have
\begin{equation}
\label{Est}
\sum_{kb\le m\le n}\mu(k)\langle bt \rangle=O(n^2).
\end{equation}
Note that the implied constant does not depend on $t$.
From (\ref{B11}),(\ref{T1_0}) and (\ref{Est}), we have shown
\begin{equation}
\label{T1}
B(n,t,1)=\frac{tn^3}{\pi^2}+ O(n^2).
\end{equation}
We do not know whether
$$
\sum_{kb\le m\le n}\mu(k)\langle bt \rangle=o(n^2)
$$
holds for all $t$. 
If this estimate is valid for a fixed $t$, we observe the second main term:
$$
B(n,t,1)=\frac{t(n^3+3n^2)}{\pi^2}+ o(n^2).
$$
Fujii \cite{Fujii:07} elaborated the improvement of (\ref{Est0})
by using Hecke's Dirichlet series \cite{Hecke:22}:
$$
Z_t(s):=
\sum_{b=1}^{\infty} \frac {\langle bt\rangle-\frac 12}{b^s}.
$$
The analytic property of $Z_t(s)$ heavily depends on the Diophantine approximation 
property of $t$ by rationals. The refinement of (\ref{Est0}) 
in the proofs of \cite[Theorem 1 and 2]{Fujii:07} imply
\begin{equation}
\label{Special}
B(n,t,1)=\frac{t(n^3+3n^2)}{\pi^2}+ O\left(n^2\exp\left(-c \left(\log n\cdot \log \log n\right)^{1/3}\right)\right)
\end{equation}
for almost all $t$, including all algebraic numbers.\footnote{The error term of (\ref{Special}) 
can be replaced by the one in (\ref{B11}), 
because we need to improve (\ref{Est}) but not necessarily (\ref{Est0}).
This may be discussed elsewhere.}
Note that the proofs are rather different between rational $t$ and algebraic
irrational $t$, and the implied $O$ constant could be very sensitive to the choice of $t$.

\section{Preliminaries}
Let $\e(x)=\exp(2\pi \I x)$. 
Let $\d$ be the natural metric on the torus
$\T=\R/\Z$. For a given interval $[\alpha,\beta]\subset [0,1)$ and a positive $J>2/(1-\beta+\alpha)$
there exists a smooth function $V_J$ of period $1$ such that
\begin{eqnarray}
\nonumber
V_J(z)=1 & z\in [\alpha,\beta]\\
\label{Fourier}
V_J(z)=0 & \d(z \bmod{\Z}, [\alpha,\beta])\ge 1/J\\
\nonumber
0\le V_J(z)\le 1 & \mathrm{ otherwise},
\end{eqnarray}
whose Fourier expansion is
$$
V_J(z)=\sum_{h\in \Z} v_J(h) \e(h z)
$$
with 
\begin{equation}
\label{FC}
|v_J(h)|\le \min\left(\frac{2}{\pi|h|}, \frac{2J}{(\pi h)^2}\right),
\end{equation}
see \cite[Chapter 1,Lemma 12]{Vinogradov:80}.
We shall use a large sieve inequality \cite[Theorem 7.2]{Iwaniec-Kowalski},\cite[Lemma 2.4]{Bombieri-Iwaniec}:
\begin{lem}
\label{Large}
For any real numbers $x_m, y_m$ with $|x_m|\le X$ and $|y_m|\le Y$ and $\alpha_m, \beta_m\in \C$, we have
$$
\left|\sum_{m}\sum_n \alpha_m \beta_n \e(x_m y_n)\right|\le 
5\sqrt{1+XY}\left(
\sum_{|x_i-x_j|<1/Y} |\alpha_i \alpha_j|
\sum_{|y_i-y_j|<1/X} |\beta_i \beta_j| \right)^{1/2}.
$$
\end{lem}
We prepare two more lemmas.

\begin{lem}
\label{Sep}
$$
\left|
\frac1{\pi}\int_{-T}^T \exp(\I\, \alpha x)\,\frac{\sin{\beta x}}x dx -\delta(\beta-\alpha)
\right|
= 
O\left(\min\left(1,\frac1{T\left|\beta-\alpha\right|}\right)\right)
$$
for $T,\alpha,\beta>0$ with
$$
\delta(y)=\begin{cases} 1 & y>0\\
                        1/2 & y=0\\
                        0 & y<0.
\end{cases}
$$
\end{lem}
\proof
This follows from the sine integral formula
$$
\frac 2{\pi}\int_{0}^T \frac{\sin x}x dx =1 + O\left(\min\left(1,\frac1T\right)\right),
$$
c.f. \cite[p.166]{Davenport}, \cite[Lemma 13.11]{Iwaniec-Kowalski}.
\qed

\begin{lem}
\label{GCDSum}
$$
\sum_{\substack{H<h\le 2H\\ M<a\le 2M}}(h,a)=\frac{6HM}{\pi^2} \log(2\min(H,M)) + O(HM)
$$
\end{lem}
\proof
Putting $d=(h,a)$, we have
\begin{align*}
\sum_{\substack{H<h\le 2H\\ M<a\le 2M}}(h,a)=&\sum_{d\le 2 \min (H,M)} d \sum_{\substack{H/d<i\le 2H/d\\ M/d<j\le 2M/d}}\, 
\sum_{k \mid (i,j)} \mu(k)\\
=&\sum_{d\le 2 \min(H,M)} d \sum_{k\le \frac{2\min (H,M)}d} \mu(k) \sum_{\substack{H/(dk)<i\le 2H/(dk)\\ M/(dk)<j\le 2M/(dk)}}1,
\end{align*}
where the last sum is
$$
\frac H{dk} \frac M{dk} + O\left(\frac H{dk}\right)+ O\left(\frac M{dk}\right)+O(1).
$$
The proof follows from $\sum_{d\le X} \log(X/d)=X+O(\log(X))$.
\qed

\section{Proof of Theorem \ref{Main}}
Our strategy is to treat $B(n,t,u)$ as a function on $t$. Then 
$B(n,t,u)$ is a non-decreasing step function having a finite number of rational discontinuities.
Every gap between the discontinuities of $B(n,s,t)$ as a function on $t$
is greater than $1/n^2$.
Depending on $n$, we will choose $t_1\le t\le t_2$ that $t_2-t_1$ is small and $t_1,t_2$ have
a suitable property to estimate error terms. If we get the same asymptotic formulas
and their implied constants of the Landau and Vinogradov symbols 
are independent of this choice of $t_i$, 
we obtain the estimate for $t$ 
by the non-decreasing property, because it is sandwiched by the same formula.
Therefore if we could perturb $t$ and get the same error term (up to negligible terms), 
then we are done. 
Hereafter we shall be cautious 
on the above implied constants, whether they can be independent of the choices of $t_i$.

At the cost of an additional error term $O(n^2)$ which takes care of 
the case $i/j=t$,  we may assume that $0<t<1$, $t$ is rational and $B(n,t,u)$ is 
continuous at $t$, i.e., constant in the neighborhood of $t$.
Recall that from the middle of \S 3, $i$ and $j$ are positive integers and
$$
B(n,t,u)=1+n+
\sum_{\substack{i,j\le m\le n,\, i/j\le t\\\langle mi/j \rangle < u,\, (i,j)=1} }1.
$$

Taking  $V_J(x)$ with respect to the interval $[0,u]$, we have
\begin{equation}
\label{Init}
\sum_{\substack{i,j\le m\le n,\, i/j\le t\\(i,j)=1} } V_J\left(\frac{mi}j\right)
=B(n,t,u) -n +  O\left(\frac {n^3}{J}\right).
\end{equation}
As (\ref{Fourier}) implies $v_J(0)=u+O(1/J)$, the main term is
\begin{eqnarray}
\nonumber
\sum_{\substack{i,j\le m\le n,\, i/j\le t\\(i,j)=1} } \left(u + O\left(\frac {1}{J}\right)\right) 
&=&u \sum_{\substack{i,j\le m\le n,\, i/j\le t\\(i,j)=1} }1 +  O\left(\frac {n^3}{J}\right)\\
\label{MainT}
&=&\frac {tun^3}{\pi^2} + O(n^2) + O\left(\frac {n^3}{J}\right)
\end{eqnarray}
from (\ref{T1}). The remainder is 
\begin{eqnarray}
\nonumber
&&\sum_{0\neq h\in \Z} v_J(h) \sum_{\substack{i,j\le m\le n,\, i/j\le t} } 
\e\left(\frac{hmi}j\right) \sum_{k \mid (i,j)} \mu(k)\\
\label{MuSum}
&=&\sum_{k\le n} \mu(k)
\sum_{0\neq h\in \Z} v_J(h) \sum_{\substack{a,b\le m/k,\, a/b\le t\\ m\le n}}
 \e\left(\frac{hma}b\right)
\end{eqnarray}
with positive integers $a, b$. From (\ref{FC}), our target is to estimate
\begin{eqnarray}
\nonumber
C(n,k,t,u)&:=&\sum_{0\neq h\in \Z}v_J(h) \sum_{\substack{a,b\le m/k,\, a/b\le t\\m\le n}}
\e\left(\frac{hma}b\right)\\
\label{Target}
&=&
\sum_{0\neq |h|\le H}v_J(h) \sum_{\substack{a,b\le m/k,\, a/b\le t\\ m\le n}} \e\left(\frac{hma}b\right)
+O\left(\frac {Jn^3}{Hk^2}\right).
\end{eqnarray}
By Lemma \ref{Sep} and (\ref{FC}), we have
\begin{eqnarray}
&&\int_{-U}^U
\int_{-T}^T\sum_{1\le h\le H}v_J(h) \sum_{\substack{a,b\le n/k\\ m\le n}} 
\e\left(\frac{hma}b\right) \left(\frac ab\right)^{\I x} \left(\frac bm\right)^{\I y} 
\label{Int0}
\\
\nonumber
&&\hspace{5cm} \times \frac {\sin(x \log(t))}{\pi x} \frac {\sin(y \log(1/k))}{\pi y} dx dy
\\
\label{Int}
&=&
\sum_{1\le h\le H}v_J(h) \sum_{\substack{a,b\le m/k,\, a/b\le t\\m\le n}} \e\left(\frac{h ma}b\right)\\
&&
\nonumber
+O\left(\log U \log H \sum_{\substack{a,b\le n/k\\m\le n}}  \frac 1{T|\log \frac {tb}a|}\right)
+O\left(\log H \sum_{\substack{a,b\le n/k,\, a/b\le t\\
m\le n}}  \min\left(1,\frac 1{U|\log \frac {bk}m|}\right)\right).
\end{eqnarray}
It is important to note that the variables $h,a,b,m$ are independent 
in the sum in the double integral (\ref{Int0}).
We may further assume that $t$ has a denominator $p$, which is the smallest prime exceeding $n^2$. Then
we have $|\log (1+(\frac{tb}a -1))|\ge 1/(2np)$ and $\frac {1}{T|\log \frac{tb}a|} \le \frac{4n^3}T$
which implies
\begin{equation}
\label{IntErr}
O\left(\log U \log H \sum_{a,b\le n/k,\, m\le n} \frac 1{T|\log \frac {tb}a|}\right)
=O\left(\frac{n^6}{k^{2} T}\log U \log H \right).
\end{equation}
Since one can show that the case $bk=m$ is negligible, we have
\begin{equation}
\label{IntErr2}
O\left(\log H \sum_{\substack{a,b\le n/k,\, a/b\le t\\
m\le n}}  \min\left(1,\frac 1{U|\log \frac {bk}m|}\right)\right)
=O\left(\frac{n^4}{k^{2} U} \log H \right).
\end{equation}
To deal with (\ref{Int0}), let us apply Lemma \ref{Large} to 
$$
G(K,L,M,N):=\sum_{\substack{K<m\le 2K\\L<h\le 2L}} \sum_{\substack{M<a\le 2M\\N<b\le 2N}}
v_J(h) a^{\I x} b^{\I(y-x)} m^{-\I y}
\e\left(\frac{hma}b\right)
$$
where $M,N\le n/(2k)$, $K\le n/2$ and $K,L,M,N\in \N$. 
Define $x_m$ and $y_n$ by rearranging multi-sets 
$$\{ha\ |\ L< h\le 2L,\, M< a\le 2M\}$$
and
$$\left\{\left.\frac mb\ \right|\ K<m\le 2K,\, N< b\le 2N\right\}$$
in the non-decreasing order, keeping multiplicity. Since $|x_m|\le 4LM$, $|y_n|\le 2K/N$,
we have a bound:
\begin{align*}
&|G(K,L,M,N)|^2\\
&\le  25 \left(1+\frac{8LMK}N\right) \left(\sum_{|h_1a_1-h_2a_2|< \frac{N}{2K}} |v_J(h_1) v_J(h_2)|\right)
\left(\sum_{|m/b_1-m/b_2|<\frac 1{4LM}} 1\right) \\ 
&\ll  \left(1+\frac{8LMK}N\right) \frac 1{L^2} \left(\frac N{2K}+1\right)
\left(\sum_{\substack{L<h_1\le 2L\\ M<a_2\le 2M}} (h_1, a_2)\right) KN \left( \frac{N^2}{4LMK} +1\right)\\
&\ll \left(1+\frac{LMK}N\right) \frac{MNK}L \left(\frac N{K}+1\right) \left( \frac{N^2}{LMK} +1\right) \log \min(L,M) 
\end{align*}
by Lemma \ref{GCDSum}. 
If
$L\le N/(8MK)$, then
\begin{equation}
\label{VerySmallL}
|G(K,L,M,N)|^2 \ll \left(\frac {N^4}{KL^2}+\frac{N^3}{L^2}\right) \log \min(L,M)
\ll \frac{N^4}{L^2}\log \min(L,M).
\end{equation}
If $N/(8MK)<L\le N^2/(4MK)$, we have
\begin{equation}
\label{SmallL}
|G(K,L,M,N)|^2 \ll \left(\frac {MN^3}{L}+\frac{MN^2K}{L} \right)\log \min(L,M).
\end{equation}
When $L> N^2/(4MK)$, we get the estimate
\begin{equation}
\label{LargeL}
|G(K,L,M,N)|^2 \ll (M^2NK + M^2K^2) \log \min(L,M).
\end{equation}
Now (\ref{VerySmallL}),(\ref{SmallL}),(\ref{LargeL})
shows 
\begin{equation}
\label{Summary}
G(K,L,M,N)\ll \frac{n^2(\log n)^{1/2}}k
\end{equation}
for any $K,L,M,N$.
Set $J=\max(n,3/(1-t))$ and $H=J^2$.
Summing up 
(\ref{Target}),(\ref{Int}),(\ref{IntErr}),(\ref{IntErr2}),(\ref{Summary})
 and a similar computation for negative $h$, we obtain
$$
C(n,k,t,u)\ll \frac {n^2}{k^2} 
+ \frac {n^6}{k^2 T}(\log U)(\log n) + \frac {n^4}{k^2 U}\log n + \frac {n^2}{k} (\log T)(\log U) (\log n)^{9/2}
$$
where the increase of the exponent of the last $\log n$ comes from
the number of dissections 
$$\sum_{1\le b\le N}=\sum_{j\ge 1}\, \sum_{2^{-j}N<b\le 2^{-j+1}N}$$
and the similar ones for $a,m$ and $h$.
Taking $U=n^2$ and $T=n^4$,
we obtain
$$
C(n,k,t,u)= O\left(\frac{n^{2}}k (\log n)^{13/2}\right).
$$
Therefore from (\ref{Init}),(\ref{MainT}) and (\ref{MuSum}), we have
$$
B(n,t,u)=\frac{tun^3}{\pi^2}+ O\left(n^{2} (\log n)^{15/2}\right).
$$

\begin{rem}
\label{A}
One can apply the same method to $A(m,t,u)$. A slightly easier computation
gives
$$
A(m,t,u)= \frac{3tum^2}{\pi^2}+ O\left(m^{3/2} (\log m)^{9/2}\right)
$$
which falls short in showing the above bound for $B(n,t,u)$.
\end{rem}

{\bf Acknowledgment.}
I am particularly grateful to my colleague Hiroshi Mikawa for showing the idea to 
introduce the large sieve inequality in this setting. Thanks are also due to Shin'ichi Yasutomi and Fujii Akio for 
inspiring discussion.

\end{document}